# An estimate of asymptotics of the moments of additive arithmetic functions with a limit distribution defined on a subset of the natural series

Victor Volfson

ABSTRACT We study the asymptotics of the moments of arithmetic functions that have a limit distribution, not necessarily normal, defined on a subset of the natural series that satisfies certain requirements. Several assertions are proved on estimating the asymptotics of the moments of strongly additive arithmetic functions and also with additive functions of the class H that have a limit distribution and are defined on a subset of the natural series.

## 1. INTRODUCTION

The asymptotic behavior of the moments of arithmetic functions on a subset of the natural series is related to the distribution of sums of functions of prime numbers on a subset of the natural series [1].

The asymptotic behavior of the moments of arithmetic functions on an arithmetic progression was considered in [2].

The paper investigates the asymptotics of the moments of arithmetic functions with limit distribution on a subset of the natural series that satisfies certain requirements.

Let's start with basic definitions. An integer sequence is a strictly increasing sequence of natural numbers. An example of an integer sequence is the sequence of prime numbers: 2, 3, 5, 7,…

Based on the asymptotic law of primes [3], the density of primes is determined by the formula:

$$d(x) = \pi(x)/x = |\{n \leq x : n \in P\}|/x = \frac{1}{\ln x}(1 + o(1)), \qquad (1.1)$$

where $P$ is the set of primes, $\pi(x)$ is the number of primes not exceeding $x$.

Let's consider a subset of the natural series $g$, which is also an integer sequence.

Let us denote the asymptotic density of natural numbers belonging to the intersection $g \cap P$:

$$d^*(g) = \lim_{x \to \infty} \pi(g, x)/\pi(x) = \lim_{x \to \infty} |\{n \leq x : n \in g \cap P\}|/|\{n \leq x : n \in P\}|. \qquad (1.2)$$

If $d^*(g)$ is not equal to zero, then having in mind (1.2) we obtain:

$$\lim_{x \to \infty} \frac{\pi(g, x)}{d^*(g)\pi(x)} = 1. \qquad (1.3)$$

Based on (1.1) and (1.3), the asymptotics of the number of natural numbers that do not exceed the value $x$, belonging to $g \cap P$, is equal to:

$$\pi(g,x) = d^*(g)\pi(x)(1+o(1)) = d^*(g)\frac{x(1+o(1))}{\ln x}. \tag{1.4}$$

Therefore, determining the asymptotics of the number of natural numbers that do not exceed the value $x$, belonging to $g \cap P$, in the case when $d^*(g)$ is not equal to zero, reduces to finding the asymptotic density $d^*(g)$ using formula (1.2).

For example, based on [4] and (1.2), the asymptotic density primes on a subset of the natural series $g$, which is a polynomial of degree $h$, is determined by the formula:

$$d^*(g) = \lim_{x\to\infty} \pi(g,x)/\pi(x) = \lim_{x\to\infty} \frac{c(g)}{h}\int_2^x \frac{dt}{\ln t} / \int_2^x \frac{dt}{\ln t} = \frac{c(g)}{h},$$

where $c(g) = \prod_p \frac{1-\omega(p)/p}{1-1/p}$, $\omega(p)$ is the number of comparison solutions $g \equiv 0 \pmod{p}$.

Therefore, in this case, as $\omega(p) < p$, then $d^*(g)$ is not equal to zero, and therefore satisfies the specified requirements.

If a subset of the natural series is an arithmetic progression $g = kn + l\,(k,l=1)$, then the value of the asymptotic density $d^*(g) = 1/\varphi(k)$ (where $\varphi(k)$ is the value of the Euler function at the point $k$) is not equal to zero.

An arithmetic function is a function defined on a subset of natural numbers and taking values on the set of complex numbers.

An arithmetic function is additive if it satisfies the condition:
$$f(m) = f(p_1^{a_1}...p_t^{a_t}) = f(p_1^{a_1}) + ... + f(p_t^{a_t}) = \sum_{p^\alpha \| m} f(p^\alpha).$$

Problems arise in determining an asymptotic of the mean value of arithmetic functions [5], and even more so in determining the asymptotics of moments of higher orders. Therefore, I propose to approach the solution of this problem from the other side.

The premises of this approach are as follows.

1. Let there are two sequences of random variables that have the same limit distribution. Naturally, these sequences have the same asymptotics of the mean value and moments of higher orders.

2. It is known that an arithmetic function can be represented as a sequence of random variables. If this sequence converges to a limit distribution function that coincides with the limit distribution function of another sequence of random variables, then the asymptotics of the mean value and the moments of higher orders coincide for the arithmetic function and the sequence of random variables. It is possible to construct a sequence of random variables in this case, in which the indicated characteristics are determined more simply. On the other hand, the given sequence of random variables must converge to the same distribution function as the arithmetic function. Then one can more simply determine the asymptotics of the characteristics of an arithmetic function from the asymptotics of the characteristics of a given sequence of random variables.

3. If two arithmetic functions have the same limit distribution, then they have the same asymptotics of all moments.

Let's start with strongly additive arithmetic functions. Let me remind you that, by definition, a strongly additive arithmetic function is a function for which $f(p^a) = f(p)$. Therefore, for an arbitrary natural number $m = p_1^{a_1}...p_t^{a_t}$ for a strongly additive arithmetic function, we have:

$$f(m) = f(p_1^{a_1}...p_t^{a_t}) = f(p_1^{a_1}) + ... + f(p_t^{a_t}) = f(p_1) + ... + f(p_t) = \sum_{p|m} f(p). \qquad (1.5)$$

Any initial segment of the natural series $\{1,...,n\}$ can be naturally transformed into a probability space $(\Omega_n, \mathcal{A}_n, \mathbb{P}_n)$ by taking $\Omega_n = \{1,...,n\}$, $\mathcal{A}_n$ — all subsets of $\Omega_n$, $\mathbb{P}_n(A) = \frac{\#A}{n}$. Then an arbitrary (real) function $f$ of the natural argument $m$ (more precisely, its restriction to $\Omega_n$) can be considered as a random variable $\xi_n$ on this probability space: $\xi_n(m) = f(m)$, $1 \leqslant m \leqslant n$.

Therefore, we can write the Chebyshev inequality for an arithmetic function $f(m), m = 1,...,n$ on a given probability space:

$$P_n(|f(m) - A_n| \leqslant b\sigma_n) \geqslant 1 - 1/b^2. \qquad (1.6)$$

where the value $b \geqslant 1$ and $A_n, \sigma_n$, respectively, the mean value and the standard deviation of $f(m), m = 1,...,n$.

Let us put $b = b(n)$ in (1.6), where $b(n)$ is an indefinitely increasing function for the value $n \to \infty$:

$$P_n(|f(m) - A_n| \leq b(n)\sigma_n) \geq 1 - 1/b^2(n). \tag{1.7}$$

The limit $P_n$ in (1.7) exists when the value $n \to \infty$:

$$P_n(|f(m) - A_n| \leq b(n)\sigma_n) \to 1, n \to \infty. \tag{1.8}$$

Expression (1.8) is an analogue of the law of large numbers for arithmetic functions [6].

Turan proved [7] that if an arithmetic function $f(m), m = 1,...,n$ is strongly additive, satisfies the condition for all primes $p$: $0 \leq f(p) < c$ and $A_n \to \infty$ for the value $n \to \infty$, then the following form of the analog of the law of large numbers holds for it:

$$P_n(|f(m) - A_n| \leq b(n)\sqrt{A_n}) \to 1, n \to \infty. \tag{1.9}$$

So, each arithmetic function $f(m), m = 1,...,n$ can be associated with a sequence of random variables $\xi_n$ (living in different probability spaces). Each such random variable has a distribution function, so it is natural to consider the question of the convergence of the sequence of distribution functions corresponding to the arithmetic function to the limit distribution function.

In 1939, Erdős and Katz proved [8] that the limit distribution function for an arithmetic function $\omega(m), m \leq n$ is normal, i.e., for any $x \in R$ one runs:

$$P_n\left(\frac{\omega(m) - \ln\ln n}{\sqrt{\ln\ln n}} \leq x\right) \to \Phi(x), n \to \infty. \tag{1.10}$$

where $\Phi(x)$ is the standard normal distribution function.

Erdős and Katz also generalized (1.10) and proved that if $f(m)$ is a strongly additive arithmetic function and it holds $|f(p)| \leq 1$ for all primes $p$, then for $\sigma_n \to \infty, n \to \infty$:

$$P_n\left(\frac{f(m) - A_n}{\sigma_n} \leq x\right) \to \Phi(x), n \to \infty. \tag{1.11}$$

It is also interesting to consider additive arithmetic functions that have a limit distribution other than normal.

We define a random variable $f^{(p)}(m) = f(p)$ for each prime $p(p \leq n)$, $f^{(p)}(m) = f(p)$ if $p \mid m$ and otherwise $f^{(p)}(m) = 0$.

Then $f^{(p)}(m) = f(p)$ with probability $\frac{1}{n}[\frac{n}{p}]$ and $f^{(p)}(m) = 0$ with probability $1 - \frac{1}{n}[\frac{n}{p}]$.

Therefore, the average value $f^{(p)}(m)$ over the interval $[1, n]$ is:

$$E[f^{(p)}, n] = \frac{f(p)}{n}[\frac{n}{p}]$$

It is performed for strongly additive $f(m)$:

$$f(m) = \sum_{p \leq n} f^{(p)}(m). \tag{1.12}$$

Having in mind (1.12), the average value $f(m)$ on the interval $[1, n]$ is:

$$A_n = \sum_{p \leq n} \frac{f(p)}{n}[\frac{n}{p}].$$

Therefore, the following asymptotic of the mean value of a strongly additive arithmetic function $f(m)$ holds:

$$A_n \to \sum_{p \leq n} \frac{f(p)}{p}, n \to \infty. \tag{1.13}$$

It was shown [6] that the asymptotic of the dispersion of a strongly additive arithmetic function $f(m)$ is determined by the formula:

$$D(n) \sim \sum_{p \leq n} \frac{f^2(p)}{p}. \tag{1.14}$$

Based on (1.12), any strongly additive arithmetic function $f(m), m = 1, ..., n$ can be represented as a sum of random variables:

$$f(m) = \sum_{p \leq n} f^{(p)}(m).$$

The question arises, in which case will the random variables $f^{(p)}(m)$ be asymptotically independent? This is done for the so-called "stripped" arithmetic functions $f(m)_r$ [6].

Let the function $r = r(n)$ grow more slowly than any positive power $n$ for $n \to \infty$. Denote the mean value and variance for the function $r = r(n)$, respectively: $A(r), D(r)$.

Then it was proved that for a real additive arithmetic function $f(m), m = 1,...,n$ with mean value and variance, respectively $A(n), D(n)$, there is such an indefinitely increasing function $r = r(n)$, provided that $\ln r(n) / \ln n \to 0$, $D(r)/D(n) \to 1$ and limit distribution laws:

$$P_n\{\frac{f(m)_r - A(r)}{D(r)} < x\}, P_n\{\frac{f(m) - A(n)}{D(n)} < x\}$$

exist only simultaneously and coincide in this case. We assume that $f(m)$ belongs to the class $H$ of additive arithmetic functions in this case.

It is shown [6] that strongly additive arithmetic functions $f^*(m)$ defined on a subset of the natural series g can act as additive arithmetic functions $f(m)$ of the class $H$ in the case when $d^*(g)$ is not equal to zero.

It is also proved that a strongly additive arithmetic function $f(m)$ belongs to the class $H$ if its variance satisfies the condition: $\ln D(n) = o(\ln \ln(n))$ and the following theorems are proved.

Theorem 1.1

In order for the distribution laws $P_n\{\frac{f(m) - A(n)}{\sqrt{D(n)}} < x\}$, where $f(m)$ is a strongly additive arithmetic function of the class $H$, to converge to the limit with a variance equal to 1, it is necessary and sufficient that there is such a no decreasing function $K(u)$ with a variation equal to 1 that for at $n \to \infty$ all points of continuity $K(u)$ the following condition is satisfied:

$$\frac{1}{D(n)} \sum_{p \leq n, f(p) < u\sqrt{D(n)}} \frac{f^2(p)}{p} \to K(u),$$

where, respectively, $A(n), D(n)$ are the mean value and variance $f(m)$, and $K(u)$ is the Kolmogorov function:

$$K(u) = 0, if u < A; K(u) = \mu(1 - u^2 / A^2), if A \leq u < 0;$$

$$K(u) = vu^2 / C^2 + 1 - v, if 0 < u \leq C; K(u) = 1, if u > C,$$

where $A, C, \mu, v$ are constants and $\mu \geq 0, v \geq 0, \mu + v \leq 1, A \leq 0, C \geq 0$.

Theorem 1.2

The class of limit laws to which the distribution law $P_n\{\frac{f(m) - A(n)}{\sqrt{D(n)}} < x\}$ tends with variance equal to 1 for a strongly additive arithmetic function $f(m)$ belonging to the class $H$ coincides with the class of laws $K(u)$.

It was said above that if there is a limit distribution of a strongly additive arithmetic function, then it coincides with the limit distribution of the sum of asymptotically independent random variables, and this arithmetic function belongs to the class $H$ in this case.

It is known [9] that the limit distribution of the sum of asymptotically independent random variables that take two values, under certain conditions, which we will consider in detail below, coincides with one of the distribution laws $K(u)$.

Therefore, if there is a limit distribution of a strongly additive arithmetic function, then it coincides with one of the distribution laws $K(u)$.

2. ESTIMATE OF ASYMPTOTICS OF THE MOMENTS OF ARITHMETIC FUNCTIONS WITH LIMIT DISTRIBUTION, DEFINED ON A SUBSET OF THE NATURAL SERIES

Assertion 2.1

The distribution laws $P_n\{\frac{f(m) - A(n)}{\sqrt{D(n)}} < x\}$, where $f(m)$ is a strongly additive arithmetic function of the class $H$ defined on a subset of the natural series $g$, converge to the limit function with a variance equal to 1, if there exists such a non-decreasing function $K(u)$ with a variation equal to 1, so that at all points of continuity $K(u)$ the following condition is satisfied:

$$\frac{1}{D(n)} \sum_{p \leq n, p \in g, f(p) < u\sqrt{D(n)}} \frac{f^2(p)}{p} \to K(u),$$

where $A(n), D(n)$, respectively, are the mean value and variance of $f(m)$, and $K(u)$ is the Kolmogorov function:

$$K(u) = 0, \text{if } u < A; K(u) = \mu(1 - u^2/A^2), \text{if } A \leq u < 0;$$

$$K(u) = vu^2/C^2 + 1 - v, \text{if } 0 < u \leq C; K(u) = 1, \text{if } u > C,$$

where $A, C, \mu, v$ are constants and $\mu \geq 0, v \geq 0, \mu + v \leq 1, A \leq 0, C \geq 0$.

In addition, the asymptotic density $d^*(g)$ must be non-zero. This condition is necessary and sufficient.

Proof

Based on the introduction, strongly additive arithmetic functions $f^*(m)$ defined on a subset of the natural series $g$ can act as additive arithmetic functions $f(m)$ of the class $H$ in the case when $d^*(g)$ is not equal to zero.

Therefore, based on Theorem 1.1, this assertion holds.

Assertion 2.2

The class of limit laws to which the distribution law $P_n\{\frac{f(m) - A(n)}{\sqrt{D(n)}} < x\}$ tends with variance equal to 1 for a strongly additive arithmetic function $f^*(m)$ defined on a subset of the natural series $g$ (in the case when the asymptotic density $d^*(g)$ is not equal to zero) and belonging to the class $H$ coincides with the class of laws $K(u)$.

Proof

Based on the introduction, strongly additive arithmetic functions $f^*(m)$ defined on a subset of the natural series $g$, can act as additive arithmetic functions $f(m)$ of the class $H$ in the case when $d^*(g)$ is not equal to zero.

Therefore, having in mind Theorem 1.2, this assertion holds.

Assertion 2.3

Let there is a strongly additive arithmetic function - $f(m)$, defined on a subset of the natural series $g$ (in the case when the asymptotic density $d*(g)$ is not equal to zero), for which the series $\sum_{p \leq n, p \in g} \frac{f^2(p)}{p}$ - diverges. Let us denote $B(n) = (\sum_{p \leq n, p \in g} \frac{f^2(p)}{p})^{1/2}$, then, if $\frac{1}{B^2(n)} \sum_{p \leq n, p \in g} \frac{f(p)}{p} \to K(u)$ when the value $n \to \infty$, where $K(u)$ is the Kolmogorov function and $P\{\frac{|f(p)|}{B(n)} > \epsilon\} \to 0, n \to \infty$, then:

1. It is possible to construct a random variable that is the sum of independent random variables taking two values, which has the same limit distribution as $f(m)$, for $m \in g(n)$ and $n \to \infty$, for which the asymptotics of the mean value and variance, respectively, are:

$$A(n) \sim \sum_{p \leq n, p \in g} \frac{f(p)}{p}, D(n) \sim \sum_{p \leq n, p \in g} \frac{f^2(p)}{p}.$$

2. The asymptotics of the central moment of the $u$-th order for $f(m)$, at value $m \in g(n)$ and $n \to \infty$ is determined by the formula: $\sum_{p \leq n, p \in g} \frac{f^u(p)}{p}$.

Proof

Taking into account that the conditions of Assertion 2.1 are satisfied ($\frac{1}{B^2(n)} \sum_{p \leq n, p \in g} \frac{f(p)}{p} \to K(u)$, where $K(u)$ is the Kolmogorov function), then, based on this theorem, the distribution laws $P_n\{\frac{f(m) - A(n)}{\sqrt{D(n)}} < x\}$, where $f(m)$ is a strongly additive arithmetic function of the class $H$ at $m \in g(n)$ and $n \to \infty$, converge to the limit distribution with a variance equal to 1. Based on Assertion 2.2, the class of these limit laws matches with $K(u)$.

On the other hand, since the series $\sum_{p \leq n, p \in g} \frac{f^u(p)}{p}$. - diverges and $B(n) \to \infty$, then, based on Lemma 4.6 [9], there is a sequence of independent random variables that take no more than two values - $X_p$, for which the distribution law $\frac{1}{B(n)} \sum_{p \leq n, p \in g} X_p - A(n)$, under the condition $P\{\frac{|X_p|}{B(n)} > \epsilon\} \to 0, n \to \infty$, converges to the limit distribution coinciding with $K(u)$.

Let's build such a sequence of independent random variables $X_p$ and show that the condition $P\{\frac{|X_p|}{B(n)} > \epsilon\} \to 0, n \to \infty$ will be satisfied for it.

Let a random variable $X_p$ take two values: $X_p = f(p)$ with probability $1/p$ and $X_p = 0$ with probability $1 - 1/p$.

Then the mean and variance $X_p$ are respectively equal:

$$E[X_p] = f(p)/p, D[X_p] = f^2(p)/p .$$

Taking into account that by the condition $P\{\frac{|f(p)|}{B(n)} > \epsilon\} \to 0, n \to \infty$, and having in mind that $|X_p| \leq |f(p)|$, we get $P\{\frac{|X_p|}{B(n)} > \epsilon\} \to 0, n \to \infty$.

Let us define a random variable $S_n = \sum_{p \leq n, p \in g} X_p$, where $X_p$ are independent random variables. Then the mean and variance $S_n$ are respectively equal:

$$E[S_n] = \sum_{p \leq n, p \in g} \frac{f(p)}{p}, D[S_n] = \sum_{p \leq n, p \in g} \frac{f^2(p)}{p},$$

those are equal to the asymptotics, respectively, of the mean value and variance $f(m)$ at $m \in g(n)$ and $n \to \infty$. (see 1.13 and 1.14). Thus, the limit distributions $f(m)$ and $S_n$ coincide for $m \in g(n)$ and $n \to \infty$.

Therefore, we have proved the first part of the assertion. Now let's prove the second part. Since the limit distributions $f(m)$ and $S_n$ for $m \in g(n)$ and $n \to \infty$ coincide, the central moments of all orders also coincide. Therefore, it suffices to determine the central moments of the random variable $S_n$.

First, we determine the central moment of the $u$-th order for $X_p$:

$$E[(X_p - \frac{f(p)}{p})^u] = E[X_p]^u - uE[X_p^{u-1}]\frac{f(p)}{p} + \frac{u(u-1)}{2}E[X_p^{u-2}]\frac{f^2(p)}{p} - \ldots + (-1)^u \frac{f^u(p)}{p} = \frac{f^u(p)}{p} - k\frac{f^{u-1}(p)}{p}\frac{f(p)}{p} + \frac{u(u-1)}{2}\frac{f^{u-2}(p)}{p}\frac{f^2(p)}{p^2} + \ldots (-1)^u \frac{f^u(p)}{p^u}$$

Thus, the central moment of the $u$-th order for $X_p$ is equal to:

$$E[(X_p - \frac{f(p)}{p})^u] = \frac{f^u(p)}{p} + O(\frac{f^u(p)}{p^2}).$$

We get that the central moment of the $u$-th order for $S_n$ is equal to:

$$E[S_n^u, n] = \sum_{p \leq n, p \in g} \frac{f^u(p)}{p} + O(\sum_{p \leq n, p \in g} \frac{f^u(p)}{p^2}) = \sum_{p \leq n, p \in g} \frac{f^u(p)}{p} + o(\sum_{p \leq n, p \in g} \frac{f^u(p)}{p}),$$

those $E[S_n^u, n] \sim \sum_{p \leq n, p \in g} \frac{f^u(p)}{p}$.

Therefore, we have proved the second part of the assertion.

Let us consider the cases when the condition $P(\frac{|f(p)|}{B(n)} > \epsilon) \to 0, n \to \infty$ is satisfied.

If $|f(p)| \leq C$, then using the results of [1] (in the case when the asymptotic density $d^*(g)$ is not equal to zero) we obtain:

$$B(n) \sim (\sum_{p \leq n, p \in g} \frac{f^2(p)}{p})^{1/2} \leq (\sum_{p \leq n, p \in g} \frac{C^2}{p})^{1/2} = d^*(g)(\ln\ln(n))^{1/2},$$

therefore, in this case, the condition $P(\frac{|f(p)|}{B(n)} > \epsilon) \to 0, n \to \infty$ is satisfied.

If $\sup_{p \leq n, p \in g} |f(p)| = o(B(n))$, then the condition $P(\frac{|f(p)|}{B(n)} > \epsilon) \to 0, n \to \infty$ is also satisfied.

Let me remind you that in these cases the strongly additive arithmetic function $f(m)$ for $m \in g(n)$ and $n \to \infty$ has a limiting normal distribution.

We present a strongly additive arithmetic function as a more general case that, as shown in [6], has a limit distribution function $K(u)$, when the asymptotic density $d*(g)$ is not equal to zero:

$$f(p) = \sqrt{2(1+\mu sgn(A) - \nu sgn(C))\ln\ln(p)}, \text{ if } p \in Q_0;$$

$$f(p) = A\ln\ln(p), \text{ if } p \in Q_1;$$

$$f(p) = C\ln\ln(p), \text{ if }, p \in Q_2$$

where $Q_0, Q_1, Q_2$ are classes of prime numbers $p \leq n, p \in g(n)$.

Using [1], the number of primes $p \in Q_1$ is asymptotically:

$$\frac{2d*(g)\mu n}{A^2 \ln(n)\ln\ln(n)},$$

if $A\mu$ is not equal to zero.

The number of primes $p \in Q_2$ is asymptotically equal to:

$$\frac{2d*(g)\nu n}{C^2 \ln(n)\ln\ln(n)},$$

if $C\nu$ is not equal to zero.

If $A\mu = 0$ or $C\nu = 0$, then $Q_1$ and $Q_2$ are empty.

The number of primes $p \in Q_0$ is asymptotically equal to:

$$\frac{d*(g)n}{\ln(n)}.$$

Using the results of [1], we obtain the asymptotics of the central moment of $u$-th order for a given strongly additive arithmetic function:

$$\sum_{p \leq n, p \in Q_0, p \in g} \frac{f^u(p)}{p} = d*(g)(1 + \mu sgn(A) - \nu sgn(C))(\ln\ln(n))^u(1 + o(1)),$$

$$\sum_{p \le n, p \in Q_1, p \in g} \frac{f^u(p)}{p} = -d*(g)\mu sgn(A)(\ln\ln(n))^u(1+o(1)),$$

$$\sum_{p \le n, p \in Q_2, p \in g} \frac{f^u(p)}{p} = \nu sgn(C)(\ln\ln(n))^u(1+o(1)).$$

Thus, we get:

$$\sum_{p \le n, p \in g} \frac{f^u(p)}{p} = \ln\ln(n))^u(1+o(1)).$$

Let's check the condition $P(\frac{|f(p)|}{B(n)} > \epsilon) \to 0, n \to \infty$ in this example:

$$P(\frac{|f(p)|}{B(n)} > \epsilon) = \sum_{i=0}^{2} \frac{\#\{p \in Q_i, p \le n, \frac{|f_i(p)|}{B_i(n)} > \epsilon\}}{n} = O(\frac{1}{\ln(n)\sqrt{\ln\ln(n)}}),$$

those. tends to zero when $n \to \infty$.

We have considered the case (in Assertion 2.3) when the series $\sum_{p \le n, p \in g} \frac{f^2(p)}{p}$ - diverges, where $f(m)$ - is a strongly additive arithmetic function defined on a subset of the natural series $g$ (in the case when the asymptotic density $d*(g)$ is not equal to zero). Now we consider the case when the given series converges. In this case, the following assertion is true.

Assertion 2.4

Let there is a strongly additive arithmetic function $f(m)$ defined on a subset of the natural series $g$ (in the case when the asymptotic density $d*(g)$ is not equal to zero), for which the series $\sum_{p \le n, p \in g} \frac{f^2(p)}{p}$ converges and the asymptotics of the mean value and variance for $f(m), m = g(n)$ when $n \to \infty$ are equal, respectively:

$$\sum_{p \le n, p \in g} \frac{f(p)}{p}, \sum_{p \le n, p \in g} \frac{f^2(p)}{p} \quad \text{then:}$$

1. You can construct a random variable that is the sum of independent random variables that take two values, has the same limit distribution as $f(m), m = g(n)$ for $n \to \infty$, for which the asymptotics of the mean value and variance are respectively equal:

$$A(n) \sim \sum_{p \leq n, p \in g} \frac{f(p)}{p}, D(n) \sim \sum_{p \leq n, p \in g} \frac{f^2(p)}{p}.$$

2. The asymptotics of the central moment of the $u$-th order for $f(m), m = g(n)$ at $n \to \infty$ is determined by the formula: $\sum_{p \leq n, p \in g} \frac{f^u(p)}{p}$.

Proof

It was proved [6] that for a strongly additive arithmetic function $f(m), m = g(n)$ (in the case when the asymptotic density $d^*(g)$ is not equal to zero) when the series $\sum_{p \leq n, p \in g} \frac{f^2(p)}{p}$ converges, then $f(m), m = g(n)$ at $n \to \infty$ has a limit distribution.

Based on the introduction, strongly additive arithmetic functions $f(m)$ defined on a subset of the natural series g can act as additive arithmetic functions of the class $H$ in the case when $d^*(g)$ is not equal to zero.

The limit distribution $f(m), m = g(n)$ for $n \to \infty$ coincides with the limit distribution of the sum of random variables $\sum_{p \leq n, p \in g} X_p$ for $n \to \infty$ in the case, where $X_p$ are independent random variables, each taking two values.

Let's construct this random variable. Let a random variable $X_p$ take two values: $X_p = f(p)$ with probability $1/p$ and $X_p = 0$ with probability $1 - 1/p$.

Then the mean value and variance $X_p$ are respectively equal: $E[X_p] = f(p)/p, D[X_p] = f^2(p)/p$.

Let us define a random variable $S_n = \sum_{p \leq n, p \in g} X_p$, where $X_p$ are independent random variables. Then the mean and variance $S_n$ are respectively equal:

$E[S_n] = \sum\limits_{p \leq n, p \in g} \dfrac{f(p)}{p}, D[S_n] = \sum\limits_{p \leq n, p \in g} \dfrac{f^2(p)}{p}$ i.e. are equal to the asymptotics, respectively, of the mean value and variance $f(m), m = g(n)$ at $n \to \infty$. (see 1.13 and 1.14). Thus, the limit distributions $f(m), m = g(n)$ and $S_n$ coincide.

Therefore, we have proved the first part of the assertion. Now let's prove the second part. Since the limit distributions $f(m), m = g(n)$ and $S_n$ for $n \to \infty$ coincide, the central moments of all orders also coincide. Therefore, it suffices to determine the central moments of the random variable $S_n$.

First, we determine the central moment of the $u$-th order for $X_p$:

$$E[(X_p - \dfrac{f(p)}{p})^u] = E[X_p^u] - uE[X_p^{u-1}\dfrac{f(p)}{p}] + \dfrac{u(u-1)}{2}E[X_p^{u-2}\dfrac{f^2(p)}{p}] - \ldots + (-1)^u \dfrac{f^u(p)}{p} = \dfrac{f^u(p)}{p} - u\dfrac{f^{u-1}(p)}{p}\dfrac{f(p)}{p} + \dfrac{u(u-1)}{2}\dfrac{f^{u-2}(p)}{p}\dfrac{f^2(p)}{p^2} + \ldots (-1)^u \dfrac{f^u(p)}{p^u}$$

Thus, the central moment of the $u$-th order for $X_p$ is equal to:

$$E[(X_p - \dfrac{f(p)}{p})^u] = \dfrac{f^u(p)}{p} + O(\dfrac{f^u(p)}{p^2}).$$

We get that the central moment of the $u$-th order for $S_n$:

$$E[S_n^u, n] = \sum\limits_{p \leq n, p \in g} \dfrac{f^u(p)}{p} + O(\sum\limits_{p \leq n, p \in g} \dfrac{f^u(p)}{p^2}) = \sum\limits_{p \leq n, p \in g} \dfrac{f^u(p)}{p} + o(\sum\limits_{p \leq n, p \in g} \dfrac{f^u(p)}{p}),$$

those $E[S_n^u, n] \sim \sum\limits_{p \leq n, p \in g} \dfrac{f^u(p)}{p}$.

Assertion 2.5

Let there is a strongly additive arithmetic function $f(m)$ defined on a subset of the natural series $g$ (in the case when the asymptotic density $d^*(g)$ is not equal to zero), for which the series $\sum\limits_{p \leq n, p \in g} \dfrac{f^2(p)}{p}$ converges and $|f(p)| \leq C$, where $p$ is a prime number and $C$ is a constant. Then the central moment of the $u$-th order for $f(m)$ at $u \geq 2$ is bounded.

Proof

We will carry out the proof by the method of mathematical induction.

Based on Assertion 2.3, the asymptotic of the central moment of the $u$-th order for a strongly additive arithmetic function $f(m), m = g(n)$ (in the case when the asymptotic density $d^*(g)$ is not equal to zero) at $n \to \infty$ is determined by the formula:

$$\sum_{p \leq n, p \in g} \frac{f^u(p)}{p}.$$

Based on the convergence of the series $\sum_{p \leq n, p \in g} \frac{f^2(p)}{p}$, we obtain for the value $u = 2$:

$$\sum_{p \in g} \frac{f^2(p)}{p} = O(1) \text{ is the basis of induction.}$$

Assume that the asymptotic of the $u$-th central moment $f(m), m = g(n)$ at $n \to \infty$ is equal to:

$$\sum_{p \leq n, p \in g} \frac{f^u(p)}{p} = O(1).$$

Since $|f(p)| \leq C$, then:

$$\sum_{p \leq n, p \in g} \frac{|f^{u+1}(p)|}{p} \leq C \sum_{p \leq n, p \in g} \frac{|f^u(p)|}{p} = O(1).$$

Consequently:

$$\sum_{p \leq n, p \in g} \frac{|f^{u+1}(p)|}{p} = O(1) \text{ is the induction step.}$$

Example. Based on Assertion 2.5, prove that the central moments of a strongly additive arithmetic function $f(m) = \ln \frac{\varphi(m)}{m}$ are bounded for $f(m), m = g(n)$ (in the case when the asymptotic density $d^*(g)$ is not equal to zero) and $n \to \infty$.

It is true:

$$f(p) = \ln \frac{\varphi(p)}{p} = \ln \frac{p-1}{p}.$$

Let's check first the boundedness of the function $f(p)$:

$$|\ln(1-1/p)| \leq |\ln(1-1/2)| = \ln 2.$$

Let us show now the convergence of the series $\sum_{p \in g} \frac{f^2(p)}{p}$:

$$\sum_{p \in g} \frac{f^2(p)}{p} = \sum_{p \in g} \frac{(\ln(1-1/p))^2}{p} = \sum_{p \in g} \frac{(-1/p + 1/2p^2 - 2/3!p^3 + ...)^2}{p} = \sum_{p \in g}(1/p^3 + o(1/p^3))$$

- converges.

By Assertion 2.5, this function has only bounded central moments of all orders.

Assertion 2.6

Let $f(m)$ is a real additive arithmetic function defined on a subset of the natural series $m = g(n)$ (in the case when the asymptotic density $d*(g)$ is not equal to zero) and it is performed $\ln D(n) = o(\ln \ln(n))$ for a strongly additive arithmetic function $f^*(m) = \sum_{p|m} f(p)$, where $D(n)$ is the variance $f^*(m)$. Then $f(m), f^*(m)$ have the same limit distributions and asymptotics of the moments of all orders on the subset of the natural series $m = g(n)$ at $n \to \infty$.

Proof

Since the condition $\ln D(n) = o(\ln \ln(n))$ is satisfied for a strongly additive arithmetic function $f^*(m) = \sum_{p|m} f(p)$, then it belongs to the class $H$. It means that: $f^*(p^\alpha) = f^*(p) = f(p), (\alpha = 1, 2, ...)$. Therefore, $f(m)$ t also belongs to the class $H$, just like $f^*(m)$.

Consequently, $f(m), f^*(m)$ have the same limit distributions and asymptotics of the moments of all orders at $m = g(n)$ when $n \to \infty$.

Let's look at an example. Based on Assertion 2.6, find the asymptotics of all moments of the arithmetic function $f(m) = \Omega(m) - \ln \frac{\varphi(m)}{m}$ on the subset of the natural series $m = g(n)$ as $n \to \infty$ (in the case when the asymptotic density $d*(g)$ is not equal to zero).

An arithmetic function $f(m) = \Omega(m) - \ln\frac{\varphi(m)}{m}$ is the difference of two real additive arithmetic functions, so it is a real additive arithmetic function.

Since $D(n) = \sum_{p \leq n, p \in g} \frac{f^{*2}(p)}{p} \sim d^*(g)\ln\ln p$, therefore, the condition $\ln D(n) = o(\ln\ln(n))$ is satisfied and therefore $f^*(m) = \sum_{p|m}(\Omega(p) - \ln\frac{\varphi(p)}{p})$ belongs to the class $H$.

Since the conditions of Assertion 2.6 are satisfied, then $f(m), f^*(m)$ have the same limit distributions and asymptotics of the moments of all orders on the subset of the natural series $m = g(n)$, so it suffices to find the asymptotics of the central moments of the function $f^*(m) = \sum_{p|m}(\Omega(p) - \ln\frac{\varphi(p)}{p})$ for $m = g(n)$ and $n \to \infty$.

Based on assertions 2.3, 2.4 and [1], the asymptotic of the central moment of the u-th order for a strongly additive arithmetic function $f^*(m) = \sum_{p|m}(\Omega(p) - \ln\frac{\varphi(p)}{p})$ (where $f(p) = \Omega(p) - \ln\frac{\varphi(p)}{p} = 1 - \ln(1 - 1/p)$) is defined as follows:

$$\sum_{p \leq n, p \in g} \frac{f^{*u}(p)}{p} = \sum_{p \leq n, p \in g} \frac{(1 - \ln(1 - 1/p))^u}{p} = \sum_{p \leq n, p \in g} \frac{1 - u(\ln(1 - 1/p))^{u-1} + u(u-1)/2(\ln(1 - 1/p))^{u-2} + ... + (-1)^u \ln(1 - 1/p)^u}{p}.$$

Therefore, we obtain the following asymptotic of the u-order moment for the function $f^*(m)$ for $m = g(n)$ and $n \to \infty$:

$$\sum_{p \leq n, p \in g} \frac{f^{*u}(p)}{p} = \sum_{p \leq n, p \in g} \frac{1}{p} - u \sum_{p \leq n, p \in g} \frac{1}{p^k} + u(u-1)/2 \sum_{p \leq n, p \in g} \frac{1}{p^u} + ... + (-1)^u \sum_{p \leq n, p \in g} \frac{1}{p^u} = d^*(g)\ln\ln p + O(1).$$

Based on Assertion 2.6, the moments of a real additive arithmetic function $f(m) = \Omega(m) - \ln\frac{\varphi(m)}{m}$ for $m = g(n)$ and $n \to \infty$ have a similar asymptotic behavior.

Assertion 2.7

Let $f(m), f^*(m)$ are, respectively, real additive and strongly additive functions defined on a subset of the natural series $m = g(n)$ (in the case when the asymptotic density $d^*(g)$ is not equal to zero) and $f(p) = h(p)$, where $p$ is an arbitrary prime number, and the following condition $\ln D(n) = o(\ln \ln(n))$ is satisfied for the variance of strongly additive arithmetic functions: $\ln D(n) = o(\ln \ln(n))$.

Then the arithmetic functions $f(m), h(m), f^*(m), h^*(m), (h^*(m) = \sum_{p|m} h(p))$ have the same limit distributions and asymptotics of the moments of all orders on a subset of the natural series $m = g(n)$.

Proof

Since a strongly additive arithmetic function $f^*(m) = \sum_{p|m} f(p)$ satisfies the condition: $\ln D(n) = o(\ln \ln(n))$, then $f^*(m)$ belongs to the class $H$. Having in mind that $f(p) = h(p)$, then a strongly additive arithmetic function $h^*(m) = \sum_{p|m} h(p)$ also belongs to the class $H$. Therefore, the following is performed:

$$f^*(p^\alpha) = f^*(p) = f(p) = h(p) = h^*(p) = h^*(p^\alpha), \alpha = 1, 2, \ldots$$

on a subset of the natural series $m = g(n)$.

This means that real additive and strongly additive arithmetic functions $f(m), h(m), f^*(m), h^*(m)$ belong to the class $H$, and therefore have the same limit distributions and asymptotics of the moments of all orders on a subset of the natural series $m = g(n)$ (in the case when the asymptotic density $d^*(g)$ is not equal to zero).

Corollary 2.8

All real additive arithmetic functions that coincide on a subset of the natural series $m = g(n)$ (in the case when the asymptotic density $d^*(g)$ is not equal to zero) have the same strongly additive arithmetic function on a subset of the natural series $g$. If the variance of a given strongly additive arithmetic function satisfies the condition: $\ln D(n) = o(\ln \ln(n)))$, then all indicated real additive and strongly additive arithmetic functions have the same limit distributions and asymptotics of the moments of all orders on a subset of the natural series $m = g(n)$ (in the case when the asymptotic density $d^*(g)$ is not equal to zero). This means that these arithmetic functions form equivalence classes in terms of the limit distribution and asymptotics of the moments of these functions, and these equivalence classes do not intersect.

Assertion 2.9

Let there is a real additive function $f(m)$, which takes the value $f(p) = 0$ on all primes $p$ on a subset of the natural series $m = g(n)$ (in the case when the asymptotic density $d*(g)$ is not equal to zero), then $f(m)$ has bounded moments of all orders on the subset of the natural series $g$.

Proof

Let's consider a strongly additive arithmetic function $f^*(m) = \sum_{p|m} f(p)$. The value $f^*(p) = f(p) = 0$ in this case.

Therefore, the variance $f^*(m)$ (on a subset of the natural series $m = g(n)$ when the asymptotic density $d*(g)$ is not equal to zero) is equal to:

$$D(n) \sim \sum_{p \leq n, p \in g} \frac{f^*(p)}{p} = O(1).$$

Therefore, the following holds for a strongly additive arithmetic function $f^*(m)$ on a subset of the natural series $m = g(n)$ (in the case when the asymptotic density $d*(g)$ is not equal to zero):

$$\ln D(n) = o(\ln \ln(n)),$$

therefore $f^*(m)$ belongs to the class $H$.

Based on Assertions 2.3 and 2.4, the u-th moment of a strongly additive arithmetic function is equal to $\sum_{p \leq n, p \in g} \frac{f^{*u}(p)}{p}$.

Having in mind, that $f^*(p) = 0$, we obtain that the u-th moment for a strongly additive arithmetic function is equal to $\sum_{p \leq n, p \in g} \frac{f^{*u}(p)}{p} = O(1)$ on a subset of the natural series $m = g(n)$ (in the case when the asymptotic density $d^*(g)$ is not equal to zero), i.e., it is bounded.

Based on Assertion 2.8, the arithmetic functions $f(m)$ and $f^*(m)$ have the same limit distributions and asymptotics of the moments of all orders on a subset of the natural series $m = g(n)$ (in the case when the asymptotic density $d^*(g)$ is not equal to zero). Therefore, asymptotics of the moments of all orders for an arithmetic function $f(m)$ on a subset of the natural series $m = g(n)$ (in the case when the asymptotic density $d^*(g)$ is not equal to zero) are also bounded.

There is a similar theorem 4.7 [6], in which the integrality of the real additive arithmetic function is additionally required, since the proof is carried out using generating functions.

Assertion 2.9 does not require that the real additive arithmetic function be integral. For example, a real additive arithmetic function $1/2(\Omega(m) - \omega(m))$ is not integer, but satisfies all the conditions of Assertion 2.9.

Moreover, Assertion 2.9 is proved for a real arithmetic function defined on a subset of the natural series $m = g(n)$ (in the case when the asymptotic density $d^*(g)$ is not equal to zero).

Assertion 2.10

Let there is a real additive function $f(m)$ that takes the value $f(p) = 0$ on all primes $p$ on a subset of the natural series $m = g(n)$ (in the case when the asymptotic density $d^*(g)$ is not equal to zero), then the limit distribution $f(m)$ on the subset of the natural series $m = g(n)$ (in the case when the asymptotic density $d^*(g)$ is not equal to zero) coincides with the limit distribution of the real additive arithmetic function $\Omega(m) - \omega(m)$.

Proof

Let's consider a strongly additive arithmetic function $f^*(m) = \sum_{p|m} f(p)$. $f^*(p) = f(p) = 0$ in this case.

Therefore, the variance $f^*(m)$ is equal to (on a subset of the natural series $m = g(n)$ in the case when the asymptotic density is $d*(g)$ not equal to zero):

$$D(n) \sim \sum_{p \leq n, p \in g} \frac{f^*(p)}{p} = O(1) \ .$$

Therefore, the following condition is satisfied for a strongly additive arithmetic function $f^*(m)$ on a subset of the natural series $m = g(n)$ (in the case when the asymptotic density $d*(g)$ is not equal to zero):

$$\ln D(n) = o(\ln \ln(n)),$$

therefore $f^*(m)$ belongs to the class $H$.

It is performed $\Omega(p) - \omega(p) = 0$ for a real additive arithmetic function $\Omega(m) - \omega(m)$ on a subset of the natural series $m = g(n)$ (in the case when the asymptotic density $d*(g)$ is not equal to zero).

Based on Corollary 2.8, all real additive functions (in this case equal to 0) that coincide on sequences of primes on a subset of the natural series $m = g(n)$ (in the case when the asymptotic density $d*(g)$ is not equal to zero), for which the condition is satisfied for a strongly additive function $\ln D(n) = o(\ln \ln(n))$, have the same limit distribution on a subset of the natural series $m = g(n)$.

Therefore, the limit distribution of any real additive arithmetic function that takes the value 0 on a sequence of prime numbers on a subset of the natural series $m = g(n)$ (in the case when the asymptotic density $d*(g)$ is not equal to zero) coincides with the limit distribution $\Omega(m) - \omega(m)$ of a real additive arithmetic function $m = g(n)$ on a subset of the natural series $m = g(n)$ (in the case when the asymptotic density $d*(g)$ is not zero).

The limit distribution of a real additive arithmetic function $f(m) = \Omega(m) - \omega(m)$ on a natural series was studied in [10].

Assertion 2.11

Let there is a real additive arithmetic function $f(m)$ on a subset of the natural series $m = g(n)$ (in the case when the asymptotic density $d*(g)$ is not equal to zero). Then the arithmetic function $g(m) = f(m) - f^*(m)$, $f^*(m) = \sum_{p|m} f(p)$ has bounded moments of all orders on a subset of the natural series $m = g(n)$ (in the case when the asymptotic density $d*(g)$ is not equal to zero).

Proof

An arithmetic function $f^*(m)$ is a strongly additive arithmetic function on a subset of the natural series $m = g(n)$ (in the case where the asymptotic density $d*(g)$ is not equal to zero).

Since the arithmetic functions $f(m)$ and $f^*(m)$ are real additive arithmetic functions on a subset of the natural series $m = g(n)$ (in the case when the asymptotic density $d*(g)$ is not equal to zero), then $g(m) = f(m) - f^*(m)$ is also a real additive function on a subset of the natural series $m = g(n)$ (in the case when the asymptotic density $d*(g)$ is not equal to zero).

Having in mind that $f^*(p) = f(p)$, then $g(p) = f(p) - f^*(p) = 0$ for all primes $p$ on a subset of the natural series $m = g(n)$ (in the case when the asymptotic density $d*(g)$ is not equal to zero).

Thus, all the conditions of Assertion 2.9 are satisfied for a real additive arithmetic function $g(m)$, and therefore $g(m)$ has bounded moments of all orders on a subset of the natural series $m = g(n)$ (in the case when the asymptotic density $d*(g)$ is not equal to zero).

Assertion 2.12

Let $f(m)$ is a real additive arithmetic function on a subset of the natural series $m = g(n)$ (in the case when the asymptotic density $d*(g)$ is not equal to zero).

Then the limit distribution of the arithmetic function $g(m) = f(m) - f^*(m)$, $f^*(m) = \sum_{p|m} f(p)$ coincides with the limit distribution of the arithmetic function $\Omega(m) - \omega(m)$ on a subset of the natural series $m = g(n)$ (in the case when the asymptotic density $d*(g)$ is not equal to zero).

Proof

Since, respectively, $f(m), f^*(m)$ are real additive and strongly additive arithmetic functions on a subset of the natural series $m = g(n)$ (in the case when the asymptotic density $d^*(g)$ is not equal to zero), then $g(m) = f(m) - f^*(m)$ is a real additive arithmetic function on a subset of the natural series $m = g(n)$ (in the case when the asymptotic density $d^*(g)$ is not equal to zero).

Having in mind that $f^*(p) = f(p)$ then $g(p) = f(p) - f^*(p) = 0$ on a subset of the natural series $m = g(n)$ (in the case when the asymptotic density $d^*(g)$ is not equal to zero).

Thus, all the conditions of Assertion 2.11 are satisfied and the limit distribution $g(m) = f(m) - f^*(m)$ coincides with the limit distribution of an arithmetic function $\Omega(m) - \omega(m)$ on a subset of the natural series $m = g(n)$ (in the case when the asymptotic density $d^*(g)$ is not equal to zero).

Note that Assertions 2.9-2.12 hold for real additive and strongly additive arithmetic functions on a subset of the natural series $m = g(n)$ (in the case when the asymptotic density $d^*(g)$ is not equal to zero) that do not necessarily belong to the class $H$.

Let's look at examples of arithmetic functions equal to the difference between real additive and strongly additive arithmetic functions on a subset of the natural series $m = g(n)$ (in the case when the asymptotic density $d^*(g)$ is not equal to zero) that have the same limit distribution:

1. $\Omega(m) - \omega(m)$.

2. $\ln(m) - \sum_{p|m} \ln(p)$.

3. $\ln(m^s) - \sum_{p|m} \ln(p^s)$.

As we showed earlier, the additive arithmetic functions in the last two examples do not belong to the class $H$.

## 3. CONCLUSION AND SUGGESTIONS FOR FURTHER WORK

The next article will continue to study the asymptotic behavior of some arithmetic functions.

## 4. ACKNOWLEDGEMENTS

Thanks to everyone who has contributed to the discussion of this paper. I am grateful to everyone who expressed their suggestions and comments in the course of this work.